\documentclass[english]{amsart}

\usepackage[backref, colorlinks, linktocpage, citecolor = blue, linkcolor = blue]{hyperref}
\usepackage[svgnames,hyperref]{xcolor}
\usepackage{multirow}
\usepackage{enumerate}
\usepackage{lmodern} % include parenthesis in \ref
\usepackage{tikz}
\usepackage{booktabs, tabularx, multirow}
\usetikzlibrary{cd,arrows,positioning}
\tikzset{>=stealth}
\tikzcdset{arrow style=tikz}
\tikzset{link/.style={column sep=1.8cm,row sep=0.16cm}}
\tikzset{map/.style={row sep=0em, column sep=0em}}
\usetikzlibrary{shapes.geometric}
\usetikzlibrary{intersections, through}
\usepackage{verbatim}
\usepackage[normalem]{ulem}
\usepackage{amssymb}

\usepackage{array}
\usepackage{mathrsfs}
\usepackage{xspace}
\usepackage{MnSymbol}
\usepackage{mathtools}
\usepackage{stmaryrd}
\usepackage[normalem]{ulem} %\sout strike out
\usepackage{amsmath}

\usepackage{soul}
\usepackage{placeins} %\FloatBarrier

\usepackage[shortlabels]{enumitem}
\renewcommand{\sectionautorefname}{Section}

\makeatletter
\@ifdefinable\equationname{\let\equationname\equationautorefname}
\def\equationautorefname~#1\@empty\@empty\null{(#1\@empty\@empty\null)}%
\@ifdefinable\AMSname{\let\AMSname\AMSautorefname}
\def\AMSautorefname~#1\@empty\@empty\null{(#1\@empty\@empty\null)}%
\@ifdefinable\itemname{\let\itemname\itemautorefname}
\def\itemautorefname~#1\@empty\@empty\null{(#1\@empty\@empty\null)%
}%
\makeatother

\RequirePackage{amsthm}
\RequirePackage{aliascnt}
\newcommand{\basetheorem}[3]{%
	\newtheorem{#1}{#2}[#3]
	\newtheorem*{#1*}{#2}
	\expandafter\def\csname #1autorefname\endcsname{#2}
}%
\newcommand{\maketheorem}[3]{%
	\newaliascnt{#1}{#3}
	\newtheorem{#1}[#1]{#2}
	\aliascntresetthe{#1}
	\expandafter\def\csname #1autorefname\endcsname{#2}
	\newtheorem{#1*}{#2}
}%
\theoremstyle{plain}   %-------------------standard Style-------------------------

\basetheorem{theorem}{Theorem}{section}
\basetheorem{ex}{Exercise}{section}

\maketheorem{proposition}{Proposition}{theorem}
\maketheorem{fact}{Fact}{theorem}
\maketheorem{corollary}{Corollary}{theorem}
\maketheorem{lemma}{Lemma}{theorem}
\maketheorem{conjecture}{Conjecture}{theorem}
\newtheorem{definition-proposition}[theorem]{Definition-Proposition}
\theoremstyle{definition}    %------------text not italic style------------------

\maketheorem{definition}{Definition}{theorem}
\maketheorem{notation}{Notation}{theorem}
\maketheorem{defprop}{Definition-Proposition}{theorem}
\maketheorem{exercise}{Exercise}{theorem}
\maketheorem{answer}{Answer}{theorem}
\maketheorem{example}{Example}{theorem}
\maketheorem{examples}{Examples}{theorem}
\maketheorem{construction}{Construction}{theorem}

\theoremstyle{remark}    %----------------also text not italic, not bold either ------------------

\maketheorem{problem}{Problem}{theorem}

\renewcommand{\dashrightarrow}{\rat}

\newcommand{\incl}[1][r]{\ar@<-0.2pc>@{^(-}[#1] \ar@<+0.2pc>@{-}[#1]}

\title[Geometry of 3-dimensional del Pezzo fibrations in characteristic $p$]{Geometry of 3-dimensional del Pezzo fibrations in positive characteristic}

\author[F. Bernasconi]{Fabio Bernasconi}
\address{Dipartimento di Matematica Guido Castelnuovo,  Sapienza Università di Roma, Piazzale Aldo Moro, 5 I-00185 Roma}
\email{fabio.bernasconi@uniroma1.it}

\subjclass[2020]{}
\thanks{This survey was presented at the conference `Fano varieties in Cetraro' in September 2023. 
The author is partially supported by the grant PZ00P2-21610 from the Swiss National Science Foundation. }

\begin{document}

%===========================================================

\begin{abstract}
	We survey some of the recent works on the geometry of del Pezzo surfaces over imperfect fields, with applications to 3-dimensional del Pezzo fibrations in positive characteristic. 
	We place particular emphasis on cases where the general fibres of the fibrations are not smooth.
\end{abstract}

\maketitle

%===========================================================

\section{Introduction} \label{Introduction}

Let $k$ be an algebraically closed field of characteristic $p \geq 0$, and let $W$ be a smooth projective variety over $k$ whose canonical divisor $K_W$ is not pseudoeffective.
The Minimal Model Program predicts the existence of a birational modification $X$ of $W$ which is endowed with a \emph{Mori fibre space structure}, i.e. a contraction $f \colon X \to Z$ satisfying the following properties:
\begin{enumerate}
	\item $X$ is terminal $\mathbb{Q}$-factorial;
	\item $\dim(Z) < \dim(X)$ and the relative Picard rank $\rho(X/Z)=1$;
	\item $-K_X$ is $f$-ample.
\end{enumerate}

The birational geometry of $X$ can then be studied by analysing the generic fibre $X_{k(Z)}$ and the base $Z$.
While the existence of Mori fibre spaces has been proven over fields of characteristic 0 in every dimension \cite{BCHM10}, in characteristic $p>0$ only the cases of surfaces and 3-folds in characteristic $p \geq 5$ have been established \cite{Tan18, HX15, CTX15, BW17, HW22}.

In characteristic $p>0$, one has to pay additional attention to the study of generic fibres of a Mori fibre spaces as generic smoothness can fail: while the total space is smooth, all closed fibres can be singular. 
Indeed, the properties of singularities of closed fibers correspond to the properties of singularities of the \emph{geometric} generic fibers \cite[Thm 12.1.1]{Gro66}, and since the function field $k(Z)$ is in general an \emph{imperfect} field, the notions of being regular, or normal or even reduced are not stable under field extensions.
A typical example is the following:

\begin{example}
	Let $p=2$ and consider the 3-fold conic bundle:
	$$ p_u \colon X=\left\{ u_0x_0^2+u_1x_1^2+u_2x_2^2=0 \right\} \subset \mathbb{P}^2_{u} \times \mathbb{P}^2_{x} \to \mathbb{P}^2_{u}.  $$
	While $X$ is a smooth 3-fold, all closed fibres are double lines.
\end{example} 

For regular conics (which correspond to the case of the generic fibers of Mori fiber spaces of relative dimension 1), the failure of being geometrically regular occurs exclusively in characteristic 2.
In this survey, we examine the next case: the study of regular del Pezzo surfaces over imperfect fields, a topic that has received significant attention in recent years.
For us, a \emph{del Pezzo surface} over a field $K$ is a projective normal surface $X$ with $H^0(X, \mathcal{O}_X)=K$ and such that the anti-canonical divisor $-K_X$ is $\mathbb{Q}$-Cartier and ample.

\section{Bounds on the failure of geometric normality}
	
The basic tool to study geometric normality of varieties over imperfect fields is a a generalisation of Tate's base-change formula to higher dimensional varieties obtained by Patakfalvi and Waldron (with later refinements by Ji-Waldron \cite{JW21} and Tanaka \cite{Tan21}):

\begin{theorem}[{\cite[Theorem 1.1]{PW22}}] \label{cbf}
	Let $X$ be a normal variety over $K$ such that $K$ is algebraically closed in the function field $k(X)$.
	Let $f \colon Y \coloneqq (X_{\overline{k}})_{\text{red}}^{N} \to X$ be the reduced normalised base change.
	Then there exists a Weil divisor $C \geq 0$ such that $$ K_Y+(p-1)C \sim f^*K_X.$$
 \end{theorem}
	
In \cite{PW22, BT22}, the canonical bundle formula of Theorem \ref{cbf} has been applied successfully to the study of del Pezzo surfaces over imperfect fields.

\begin{theorem}[{\cite[Theorem 3.7]{BT22}}] \label{thm: geom_norm}
If $p \geq 5$, a del Pezzo surface with canonical singualarities is geometrically normal and geometrically canonical.
\end{theorem}	

The bound on the characteristic is sharp as the following example shows.

\begin{example}
Set $K=\mathbb{F}_p(s_0, s_1, s_2,s_3)$, where $p$ belongs to $\left\{2, 3\right\}$.
Consider the del Pezzo surface of degree 1 in weighted projective space:
$$X= \left\{ s_0z^2 + s_1y^3+s_2x_0^6+s_3x_1^6=0 \right\} \subset \mathbb{P}_K(1,1,2,3).$$
The surface $X$ is regular, geometrically integral but not geometrically normal. 
\end{example}

In the case of regular del Pezzo surfaces, using the classification of regular twisted forms of canonical singularities \cite{Sch08} together with the possible singularities appearing on canonical del Pezzo surfaces over algebraically closed fields one can enhance Theorem \ref{thm: geom_norm}:

\begin{corollary}
	If $p \geq 11$, a regular del Pezzo surface $X$ is geometrically regular.
\end{corollary}

Also in this case the bound on the characteristic is sharp:

\begin{example}[{\cite[Example 2.3]{MS22}}]
	Consider the weighted hypersurface 
	$$X= \left\{ w^2= tx^3y+y^3z+z^3x \right\} \subset \mathbb{P}_{\mathbb{F}_7(t)}(1,1,1,2),$$
	which is a regular del Pezzo surface of degree 2 with a geometric $A_6$-singularity.
\end{example}	
	
From the perspective of geometric non-normality, the remaining interesting cases to study arise in characteristics $p \in \left\{ 2, 3\right\}$.
For generic fibres of 3-fold Mori fibre spaces over function fields of curves, the following result of Fanelli and Schr\"oer resolves the problem completely.	
\begin{theorem}[{\cite[Theorem 14.1]{FS20a}}]
	Let $X$ a regular del Pezzo surface of Picard rank 1 over a field $K$ with $[K:K^p] \leq p$.
	Then $X$ is geometrically normal.
\end{theorem}

The following examples demonstrate that the hypotheses of the above theorem are optimal.

\begin{example}
	\begin{enumerate}
		\item The cubic surface 
		$X= \left\{ x^2y+xy^2+sz^3+tw^3=0 \right\} \subset \mathbb{P}^3_{\mathbb{F}_3(s,t)}$
		is regular but not geometrically normal.
		\item Let $Y = \left\{ x^2+sy^2+zw=0 \right\} \subset \mathbb{P}^3_{\mathbb{F}_2(s)}.$
		The blow-up at the ideal $(x^2+s, z, w)$ in the chart $(y=1)$ is a regular del Pezzo surface $X$ of degree 6 of Picard rank 2 which is not geometrically normal.
	\end{enumerate}
\end{example}

There are even more subtle examples where the irregularity $h^1(X, \mathcal{O}_X)$ of a del Pezzo surface does not vanish \cite{Sch07, Mad16}.
Maddock's example is obtained as the quotient of the (geometrically non-reduced) quadric 
$$X=\left\{ s_0X_0^2+s_1X_1^2+s_2X_2^2+s_3X_3^2=0 \right\} \subset \mathbb{P}^3_{\mathbb{F}_2(s_0, s_1, s_2, s_3)}$$ by a suitable $p$-closed foliation.
With G. Martin, we show that this is no coincidence: del Pezzo surfaces with non-vanishing irregularity have a very explicit description as quotients of foliations.

\begin{theorem}[{\cite[Proposition 4.11]{BM23}}]
	Let $X$ be a geometrically integral regular del Pezzo surface over a field $K$ of characteristic $p >0$. Assume that $h^1(X,\mathcal{O}_X) \neq 0$. 
	Then, $p=2$ $[K:K^p] \geq p^2$, $K_X^2 \leq 2$, and there exists an $\alpha_{\omega_X}$-torsor $Z \to X$ such that $Z$ satisfies the following properties:

\begin{enumerate}
	\item If $K_X^2 = 2$, then $k_Z \coloneqq h^0(Z,\mathcal{O}_Z)$ is a purely inseparable extension of $k$ of degree $2$ and $Z$ is a twisted form of $\mathbb{P}(1,1,2)$ over $k_Z$.
	\item If $K_X^2=1$, then $Z$ is a normal tame del Pezzo surface such that $\epsilon(Z/k)=1$, $K_Z^2=8$ and $(Z_{\overline{k}})_{\text{red}}^{N} \simeq \mathbb{P}^2$.
\end{enumerate}
\end{theorem}

\begin{problem}
	Does case (1) appear for regular del Pezzo surfaces?
\end{problem}

As a consequence of the above results, we prove the BAB conjecture for del Pezzo surfaces over general bases.
The geometry of klt del Pezzo surfaces has been extensively studied in a joint work with H. Tanaka in \cite{BT22}.

\begin{theorem}[\cite{Tan19}, {\cite[Theorems 5.4 and 5.12]{BM23}}]
	Fix $\varepsilon>0$ and consider the class
	$$\mathcal{X}_{{d, \varepsilon}} = \left\{X \mid X \emph{ is a } 
	\varepsilon \emph{-lc del Pezzo surface} \right\}.$$
	Then $\mathcal{X}_{\text{dP}, 1}$ (resp. $\mathcal{X}_{\text{dP}, \varepsilon}$) are bounded over $\mathbb{Z}$ (resp. $\mathbb{Z}[1/30]$).
\end{theorem}

\begin{problem}
	Are $\varepsilon$-lc del Pezzo surfaces bounded over $\mathbb{Z}$? 
	We are lacking a bound on the irregularity for klt del Pezzos in characteristic $p \leq 5$, which would generalise \cite[Theorem 1.7]{BT22}.
\end{problem}

\section{Geometry and arithmetic}

Having effectively bounded geometric non-normality to small characteristics, we turn our attention to the rationality and arithmetic properties of del Pezzo surfaces over imperfect fields.
To pursue this (and to study the 3-dimensional Cremona group over fields of positive characteristic), we develop the Sarkisov program for surfaces over imperfect fields in \cite{BFSZ24}. 
As a consequence of our classication of Sarkisov links starting from $\mathbb{P}^2_K$, we obtain the following necessary condition for rationality of del Pezzo surfaces, generalising previous results of Manin and Iskovskikh (\cite{Manin66, Isk96}).

\begin{theorem}[{\cite[Theorem 4.41]{BFSZ24}}] \label{thm: sarkisov}
	Let $X$ be a regular del Pezzo surface of Picard rank 1. 
	If $X$ is $K$-rational, then $K_X^2 \geq 5$.
\end{theorem}

In a joint work with H. Tanaka, we discuss the converse implication of Theorem \ref{thm: sarkisov}. 
The cases of geometrically non-regular del Pezzo surfaces are the most challenging, and we use explicit classical projective techniques to prove rationality.

\begin{theorem}[{\cite[Theorem 4.9]{BT24}}] \label{thm: rat_reg_dP}
	Let $X$ be a regular del Pezzo surface with $K_X^2 \geq 5$ and $X(K) \neq \emptyset$. 
	Then $X$ is $K$-rational.
\end{theorem}

The $C_1$-conjecture (see \cite{Esn23}) implies that every regular Fano variety over a $C_1$-field has a rational point. Together with Tanaka, we prove this for surfaces over imperfect fields, extending the work of Manin and Colliot-Thélène (\cite{Manin66, CT86}).

\begin{theorem}[{\cite[Theorem 4.24]{BT24}}] \label{thm: C1}
	Let $X$ be a regular del Pezzo surface over a $C_1$-field $K$. 
	Then $X(K) \neq \emptyset$.
\end{theorem}

Combining Theorems \ref{thm: rat_reg_dP} and \ref{thm: C1} we deduce the following rationality criterion for 3-dimensional del Pezzo fibrations in positive characteristic.

\begin{corollary}
	Let $k$ be an algebraically closed field and let $f \colon X \to \mathbb{P}^1_k$ be a del Pezzo fibration. 
	If $K_{X_{k(\mathbb{P}1)}} \geq 5$, then $X$ is rational.
\end{corollary}

\subsection*{Acknowledgements}

The author would like to thank A. Fanelli,  G. Martin, J. Schneider, H. Tanaka and S. Zimmermann for all the valuable discussions on the geometry of del Pezzo surfaces over imperfect fields during our collaborations.

%===========================================================

\bibliographystyle{amsalpha}
\bibliography{biblio}

\providecommand{\bysame}{\leavevmode\hbox to3em{\hrulefill}\thinspace}
\providecommand{\MR}{\relax\ifhmode\unskip\space\fi MR }
\providecommand{\MRhref}[2]{%
  \href{http://www.ams.org/mathscinet-getitem?mr=#1}{#2}
}
\providecommand{\href}[2]{#2}
\begin{thebibliography}{BCHM10}

\bibitem[BCHM10]{BCHM10}
Caucher Birkar, Paolo Cascini, Christopher~D. Hacon, and James McKernan,
  \emph{Existence of minimal models for varieties of log general type}, J.
  Amer. Math. Soc. \textbf{23} (2010), no.~2, 405--468. \MR{2601039}

\bibitem[BFSZ24]{BFSZ24}
Fabio Bernasconi, Andrea Fanelli, Julia Schneider, and Susanna Zimmermann,
  \emph{Explicit {S}arkisov program for regular surfaces over arbitrary fields
  and applications}, arXiv preprints (2024), Available at
  \href{https://arxiv.org/abs/2404.03281}{arXiv:2404.03281}.

\bibitem[BM24]{BM23}
Fabio Bernasconi and Gebhard Martin, \emph{Bounding geometrically integral del
  {P}ezzo surfaces}, Forum Math. Sigma \textbf{12} (2024), Paper No. e81, 24.
  \MR{4807862}

\bibitem[BT22]{BT22}
Fabio Bernasconi and Hiromu Tanaka, \emph{On del {P}ezzo fibrations in positive
  characteristic}, J. Inst. Math. Jussieu \textbf{21} (2022), no.~1, 197--239.
  \MR{4366337}

\bibitem[BT24]{BT24}
\bysame, \emph{Geometry and arithmetic of regular del {P}ezzo surfaces}, arXiv
  preprints (2024), Available at
  \href{https://arxiv.org/abs/2408.11378}{arXiv:2408.11378}.

\bibitem[BW17]{BW17}
Caucher Birkar and Joe Waldron, \emph{Existence of {M}ori fibre spaces for
  3-folds in {${\rm char}\,p$}}, Adv. Math. \textbf{313} (2017), 62--101.
  \MR{3649221}

\bibitem[CT87]{CT86}
Jean-Louis Colliot-Th\'{e}l\`ene, \emph{Arithm\'{e}tique des vari\'{e}t\'{e}s
  rationnelles et probl\`emes birationnels}, Proceedings of the {I}nternational
  {C}ongress of {M}athematicians, {V}ol. 1, 2 ({B}erkeley, {C}alif., 1986),
  Amer. Math. Soc., Providence, RI, 1987, pp.~641--653. \MR{934267}

\bibitem[CTX15]{CTX15}
Paolo Cascini, Hiromu Tanaka, and Chenyang Xu, \emph{On base point freeness in
  positive characteristic}, Ann. Sci. \'{E}c. Norm. Sup\'{e}r. (4) \textbf{48}
  (2015), no.~5, 1239--1272. \MR{3429479}

\bibitem[Esn22]{Esn23}
Hélène Esnault, \emph{{Rational points over C1 fields}}, arXiv preprints
  (2022), Available at
  \href{https://arxiv.org/abs/2308.08063}{arXiv:2308.08063}.

\bibitem[FS20]{FS20a}
Andrea Fanelli and Stefan Schr\"{o}er, \emph{Del {P}ezzo surfaces and {M}ori
  fiber spaces in positive characteristic}, Trans. Amer. Math. Soc.
  \textbf{373} (2020), no.~3, 1775--1843. \MR{4068282}

\bibitem[Gro66]{Gro66}
A.~Grothendieck, \emph{\'{E}l\'{e}ments de g\'{e}om\'{e}trie alg\'{e}brique.
  {IV}. \'{E}tude locale des sch\'{e}mas et des morphismes de sch\'{e}mas.
  {III}}, Inst. Hautes \'{E}tudes Sci. Publ. Math. (1966), no.~28, 255.
  \MR{217086}

\bibitem[HW22]{HW22}
Christopher Hacon and Jakub Witaszek, \emph{The minimal model program for
  threefolds in characteristic 5}, Duke Math. J. \textbf{171} (2022), no.~11,
  2193--2231. \MR{4484207}

\bibitem[HX15]{HX15}
Christopher~D. Hacon and Chenyang Xu, \emph{On the three dimensional minimal
  model program in positive characteristic}, J. Amer. Math. Soc. \textbf{28}
  (2015), no.~3, 711--744. \MR{3327534}

\bibitem[Isk96]{Isk96}
V.~A. Iskovskikh, \emph{Factorization of birational mappings of rational
  surfaces from the point of view of {M}ori theory}, Uspekhi Mat. Nauk
  \textbf{51} (1996), no.~4(310), 3--72. \MR{1422227}

\bibitem[JW21]{JW21}
Lena Ji and Joe Waldron, \emph{Structure of geometrically non-reduced
  varieties}, Trans. Amer. Math. Soc. \textbf{374} (2021), no.~12, 8333--8363.
  \MR{4337916}

\bibitem[Mad16]{Mad16}
Zachary Maddock, \emph{Regular del {P}ezzo surfaces with irregularity}, J.
  Algebraic Geom. \textbf{25} (2016), no.~3, 401--429. \MR{3493588}

\bibitem[Man66]{Manin66}
Ju.~I. Manin, \emph{Rational surfaces over perfect fields}, Inst. Hautes
  \'{E}tudes Sci. Publ. Math. (1966), no.~30, 55--113. \MR{225780}

\bibitem[MS24]{MS22}
Gebhard Martin and Claudia Stadlmayr, \emph{R{DP} del {P}ezzo surfaces with
  global vector fields in odd characteristic}, Algebr. Geom. \textbf{11}
  (2024), no.~3, 346--385. \MR{4742323}

\bibitem[PW22]{PW22}
Zsolt Patakfalvi and Joe Waldron, \emph{Singularities of general fibers and the
  {LMMP}}, Amer. J. Math. \textbf{144} (2022), no.~2, 505--540. \MR{4401510}

\bibitem[Sch07]{Sch07}
Stefan Schr\"{o}er, \emph{Weak del {P}ezzo surfaces with irregularity}, Tohoku
  Math. J. (2) \textbf{59} (2007), no.~2, 293--322. \MR{2347424}

\bibitem[Sch08]{Sch08}
\bysame, \emph{Singularities appearing on generic fibers of morphisms between
  smooth schemes}, Michigan Math. J. \textbf{56} (2008), no.~1, 55--76.
  \MR{2433656}

\bibitem[Tan18]{Tan18}
Hiromu Tanaka, \emph{Minimal model program for excellent surfaces}, Ann. Inst.
  Fourier (Grenoble) \textbf{68} (2018), no.~1, 345--376. \MR{3795482}

\bibitem[Tan21]{Tan21}
\bysame, \emph{Invariants of algebraic varieties over imperfect fields}, Tohoku
  Math. J. (2) \textbf{73} (2021), no.~4, 471--538. \MR{4355058}

\bibitem[Tan24]{Tan19}
\bysame, \emph{Boundedness of regular del {P}ezzo surfaces over imperfect
  fields}, Manuscripta Math. \textbf{174} (2024), no.~1-2, 355--379.
  \MR{4730437}

\end{thebibliography}

\end{document}